\documentclass[amstex,12pt,russian,amssymb]{article}

\usepackage{mathtext}
\usepackage[cp1251]{inputenc}
\usepackage[T2A]{fontenc}
\usepackage[russian]{babel}
\usepackage[dvips]{graphicx}
\usepackage{amsmath}
\usepackage{amssymb}
\usepackage{amsxtra}
\usepackage{latexsym}
\usepackage{ifthen}

\begin{document}

\newcounter{lemma}
\newcommand{\lemma}{\par \refstepcounter{lemma}%
{\bf Лема \arabic{lemma}.}}

\newcounter{corollary}
\newcommand{\corollary}{\par \refstepcounter{corollary}%
{\bf Наслідок \arabic{corollary}.}}

\newcounter{remark}
\newcommand{\remark}{\par \refstepcounter{remark}%
{\bf Зауваження \arabic{remark}.}}

\newcounter{theorem}
\newcommand{\theorem}{\par \refstepcounter{theorem}%
{\bf Теорема \arabic{theorem}.}}

\newcounter{proposition}
\newcommand{\proposition}{\par \refstepcounter{proposition}%
{\bf Твердження \arabic{proposition}.}}

\newcounter{example}
\newcommand{\example}{\par \refstepcounter{example}%
{\bf Приклад \arabic{example}.}}

\renewcommand{\refname}{\centerline{\bf Список літератури}}

\renewcommand{\figurename}{Мал.}

\newcommand{\proof}{{\it Доведення.\,\,}}

\medskip\medskip
{\bf В.С.~Десятка} (Житомирський державний університет імені Івана
Фран\-ка)

{\bf Є.О.~Севостьянов} (Житомирський державний університет імені
Івана Фран\-ка; Інститут прикладної математики і механіки НАН
України, м.~Слов'янськ)

\medskip\medskip\medskip
{\bf V.S.~Desyatka} (Zhytomyr Ivan Franko State University)

{\bf E.A.~Sevost'yanov} (Zhytomyr Ivan Franko State University;
Institute of Applied Ma\-the\-ma\-tics and Mechanics of NAS of
Ukraine, Slov'yans'k)

\medskip\medskip\medskip
{\bf Усувні сингулярності відображень з оберненою нерівністю
Полецького на ріманових многовидах}

{\bf Removable singularities of mappings with inverse Poletksy
inequality on Rie\-man\-nian manifolds}

\medskip\medskip\medskip\medskip
Розглядаються відкриті дискретні відображення ріманових многовидів,
які задовольняють деяку модульну нерівність. Вивчається можливість
неперервного продовження таких відображень в ізольовану точку межі.
Доведено наявність такого продовження, якщо відображення не приймає
дві та більше точок зв'язного ріманового многовиду, а мажоранта у
модульній нерівності є інтегровною по майже всіх сферах.

\medskip\medskip
We consider open discrete mappings of Riemannian manifolds that
satisfy some modulus inequality. We investigate the possibility of a
continuous extension of such mappings to an isolated point on the
boundary. It is proved that, these mappings have a specified
extension, if they omit two or more points of a connected Riemannian
manifold, and the majorant participating in the modulus inequality
is integrable over almost all spheres.

\newpage
{\bf 1. Вступ.} Як відомо, квазірегулярні відображення
$f:D\rightarrow {\Bbb R}^n,$ $n\geqslant 2,$ задовольняють умову
\begin{equation}\label{eq2}
M(\Gamma)\leqslant K N(f, D) M(f(\Gamma))
\end{equation}
для будь-якої сім'ї $\Gamma$ кривих $\gamma$ в області $D,$ де $M$
-- конформний модуль сімей кривих, $K={\rm ess \sup}\, K_O(x, f),$
$$K_{O}(x,f)\quad =\quad \left\{
\begin{array}{rr}
\frac{\Vert f^{\,\prime}(x)\Vert^n}{|J(x,f)|}, & J(x,f)\ne 0,\\
1,  &  f^{\,\prime}(x)=0, \\
\infty, & {\rm в\,інших\,випадках}
\end{array}
\right.\,,$$
$$
\Vert f^{\,\prime}(x)\Vert\,=\,\max\limits_{h\in {\Bbb R}^n
\backslash \{0\}} \frac {|f^{\,\prime}(x)h|}{|h|}\,,\quad
J(x,f)=\det f^{\,\prime}(x)\,,$$
\begin{equation}\label{eq23}
N(y, f, D)\,=\,{\rm card}\,\left\{x\in D: f(x)=y\right\}\,,
\end{equation}
$$N(f, D)\,=\,\sup\limits_{y\in{\Bbb R}^n}\,N(y, f, D)\,,$$
див., напр., \cite[теорема~3.2]{MRV$_1$}. Зокрема, квазірегулярні
відображення мають неперервне продовження в ізольовану точку межі
області за деяких умов (див., напр., \cite[теорема~2.9.3]{Ri}). У
даному рукопису ми розглянемо питання про неперервне продовження в
ізольовану точку межі більш загального класу відображень. Деякі наші
результати на цю тему були опубліковані в
\cite[теорема~5.1]{SevSkv$_1$} та \cite[теорема~1]{Sev$_2$}. Проте,
результати даної замітки в якомусь сенсі є <<найбільш загальними>>,
бо у попередніх публікаціях ми вимагали хоча б якісь додаткові
топологічні умови на відображення, крім його відкритості і
дискретності.

\medskip
Будемо вважати відомими основні об'єкти, пов'язані з рімановими
многовидами: поняття довжини та об'єму, нормального околу точки і
т.п. (див., напр., \cite{Lee}).  Усюди далі ${\Bbb M}^n$ і ${\Bbb
M}_*^n$ -- ріманові многовиди розмірності $n$ з геодезичними
відстанями $d$ і $d_*,$ відповідно,
\begin{equation}\label{eq1}
B(x_0, r)=\left\{x\in{\Bbb M}^n\,:\,d(x,x_0)<r\right\},\quad
S(x_0,r)=\left\{x\in{\Bbb M}^n\,:\, d(x,x_0)=r\right\},
\end{equation}
\begin{equation}\label{eq2B} A=A(y_0, r_1, r_2)=\{y\in
{\Bbb M}_*^n\,:\,r_1<d_*(y, y_0)<r_2\},\quad 0<r_1<r_2<r_0,
\end{equation}
$dv(x)$ и $dv_*(x)$ --- міри об'єму на ${\Bbb M}^n$ і ${\Bbb
M}_*^n,$ відповідно. Всюди далі $D$ -- область в ${\Bbb M}^n,$
$n\geqslant 2.$ Борелева функція $\rho:{\Bbb M}^n\,\rightarrow
[0,\infty] $ називається {\it допустимою} для сім'ї $\Gamma$ кривих
$\gamma$ у ${\Bbb M}^n,$ якщо
\begin{equation}\label{eq1.4}
\int\limits_{\gamma}\rho (x)\, |dx|\geqslant 1
\end{equation}
для всіх (локально спрямованих) кривих $ \gamma \in \Gamma.$ У цьому
випадку ми пишемо: $\rho \in {\rm adm} \,\Gamma .$ Нехай $p\geqslant
1,$ тоді {\it $p$-модулем} сім'ї кривих $\Gamma $ називається
величина
\begin{equation}\label{eq1.3gl0}
M_p(\Gamma)=\inf\limits_{\rho \in \,{\rm adm}\,\Gamma}
\int\limits_{{\Bbb M}^n} \rho^p (x)\,dv(x)\,.
\end{equation}
Усюди далі $M(\Gamma)=M_n(\Gamma).$ Для множин $A, B\subset {\Bbb
M}^n$ ми використовуємо позначення
\begin{equation}\label{eq1A}
{\rm dist}\,(A, B)=\inf\limits_{x\in A, y\in B}d(x, y) \,,\qquad
d(A)=\sup\limits_{x, y\in A}d(x, y)\,.
\end{equation}
Іноді замість ${\rm dist}\,(A, B)$ ми пишемо $d(A, B),$ якщо
непорозуміння неможливо.

\medskip
Нехай $x_0\in D,$ $Q\colon D\rightarrow [0,\infty]$ --- вимірна
відносно міри $v$ функція, і число $r_0>0$ є таким, що куля $B(x_0,
r_0)$ лежить у деякому нормальному околі $U$ точки $x_0$ разом із
свої замиканням. Позначимо через $S_i=S(x_0,r_i),$ $i=1,2,$
геодезичні сфери з центром у точці $x_0$ і радіусів $r_1$ і $r_2.$
Для множин $E,$ $F$ і $G$ в ${\Bbb M}^n$ позначимо через $\Gamma(E,
F, G)$ сім'ю всіх кривих $\gamma\colon[a,b]\rightarrow{\Bbb M}^n,$
які з'єднують $E$ і $F$ у $G,$ іншими словами, $\gamma(a)\in
E,\,\gamma(b)\in F$ і $\gamma(t)\in G$ при $t\in(a,\,b).$ Якщо $D$
-- область ріманового многовиду ${\Bbb M}^n,$ $f:D\rightarrow {\Bbb
M}^n$ -- деяке відображення, $y_0\in f(D)$ и
$0<r_1<r_2<d_0=\sup\limits_{y\in f(D)}d_*(y, y_0),$ то через
$\Gamma_f(y_0, r_1, r_2)$ позначимо сім'ю всіх кривих $\gamma$ в
області $D$ таких, що $f(\gamma)\in \Gamma(S(y_0, r_1), S(y_0, r_2),
A(y_0,r_1,r_2)).$ Будемо говорити, що {\it $f$ задовольняє обернену
нерівність Полецького} в точці $y_0\in f(D),$ якщо нерівність
\begin{equation}\label{eq2*A}
M(\Gamma_f(y_0, r_1, r_2))\leqslant \int\limits_{A(y_0,r_1,r_2)\cap
f(D)} Q(y)\cdot \eta^n (d_*(y,y_0))\, dv_*(y)
\end{equation}
виконується для довільної вимірної за Лебегом функції $\eta:
(r_1,r_2)\rightarrow [0,\infty ]$ такої, що
\begin{equation}\label{eqA2}
\int\limits_{r_1}^{r_2}\eta(r)\, dr\geqslant 1\,.
\end{equation}
Неважко бачити, що нерівність~(\ref{eq2*A}) принаймні у просторі
${\Bbb M}^n={\Bbb R}^n$ перетворюється у співвідношення вигляду
$M(\Gamma_f(y_0, r_1, r_2))\leqslant K\cdot M(\Gamma_f(y_0, r_1,
r_2)),$ якщо функція $Q$ обмежена числом $K\geqslant 1.$ Більше
того, якщо $f$ є відображенням з обмеженим спотворенням і має
обмежену в області $D$ функцію кратності $N(f, D),$ то ми також
маємо співвідношення~(\ref{eq2}), а отже і нерівність~(\ref{eq2*A}).
Зауважимо, що оцінки типу~(\ref{eq2*A}) добре відомі та виконуються
принаймні для $p=n$ у багатьох класах відображень (див., напр.,
\cite[теорема~3.2]{MRV$_1$}, \cite[теорема~6.7.II]{Ri} та
\cite[теорема~8.5]{MRSY}). Для $p\ne n$ схожі оцінки можна знайти,
напр., у \cite{GSU}.

Відображення $f:D\rightarrow {\Bbb M}_*^n$ називається {\it
дискретним}, якщо прообраз $\{f^{-1}\left(y\right)\}$ кожної точки
$y\,\in\,{\Bbb M}_*^n$ складається з ізольованих точок, і {\it
відкритим}, якщо образ будь-якої відкритої множини $U\subset D$ є
відкритою множиною в ${\Bbb M}_*^n.$

\medskip
У подальшому,
\begin{equation}\label{eq26}
q_{x_0}(r)=\frac{1}{r^{n-1}}\int\limits_{S(x_0,
r)}Q(x)\,d\mathcal{A}\,,
\end{equation}
де $d\mathcal{A}$ позначає елемент площі поверхні на $S(x_0, r).$
Докладно про інтеграли від функцій по поверхням ріманових многовидів
див., напр., у роботі~\cite{ARS}.

\medskip
Нехай $\overline{{\Bbb M}^n}:={\Bbb M}^n\cup\{\infty\},$ і нехай
$h:\overline{{\Bbb M}^n}\times\overline{{\Bbb M}^n}\rightarrow {\Bbb
R}$ -- метрика на $\overline{{\Bbb M}^n}.$ Припустимо, що многовид
${\Bbb M}^n$ є зв'язним, тоді $d$ -- метрика на ${\Bbb M}^n.$
Згідно \cite{SM} будемо говорити, що $h$ задовольняє умову {\it
слабкої сферикалізації,} якщо $(\overline{{\Bbb M}^n}, h)$ --
компактний метричний простір і $h$ та $d$ породжують однакову
топологію на ${\Bbb M}^n.$ Рімановий многовид ${\Bbb M}^n$
називається таким, що {\it допускає слабку сферикалізацію}, якщо
існує метрика $h:\overline{{\Bbb M}^n}\times\overline{{\Bbb
M}^n}\rightarrow {\Bbb R},$ яка задовольняє умову слабкої
сферикалізації. Аналогічно можна визначити поняття слабкої
сферикалізації для довільного метричного простору (див. \cite{SM}).

\medskip
\begin{remark}\label{rem3}
Наведене означення слабкої сферикалізації відрізняється від того, що
запропоновано в~\cite{SM}. Саме, замість умови, що $h$ та $d$
породжують однакову топологію на ${\Bbb M}^n,$ у~\cite{SM} ми
вимагаємо, щоб виконувалося наступне: $h(x, y)\leqslant d(x, y)$ при
всіх $x, y\in {\Bbb M}^n,$ де $d(x, y)$ -- геодезична метрика на
зв'язному рімановому многовиді ${\Bbb M}^n.$ При вказаному означенні
топологія на ${\Bbb M}^n,$ породжена метрикою $h,$ слабкіша
топології на ${\Bbb M}^n,$ породженою метрикою $d.$ Про рівність цих
топологій не йдеться.
\end{remark}

\medskip
Одним з найважливіших прикладів слабкої сферикалізації є {\it
хордальна (сферична) метрика} в розширеному евклідовому просторі
$\overline{{\Bbb R}^n}={\Bbb R}^n\cup\{\infty\},$
див.~\cite[Section~12]{Va}. А саме, для  $x, y\in\overline{{\Bbb
R}^n}$ покладемо $h(x, y)=|\pi(x)-\pi(y)|,$ де $\pi:\overline{{\Bbb
R}^n}\rightarrow S(e_{n+1}/2, 1/2),$ $S(e_{n+1}/2, 1/2)=\{x\in {\Bbb
R}^{n+1}:|x-e_{n+1}/2|=1/2\},$ $e_{n+1}=(0,0,\ldots,0, 1)\in {\Bbb
R}^{n+1}$ і $\pi(x)=e_{n+1}+\frac{x-{e_{n+1}}}{|x-{e_{n+1}}|^2}.$
Можна показати, що
\begin{equation}\label{eqA17}h(x,\infty)=\frac{1}{\sqrt{1+{|x|}^2}},
\ \ h(x,y)=\frac{|x-y|}{\sqrt{1+{|x|}^2} \sqrt{1+{|y|}^2}}\,, \ \
x\ne \infty\ne y\,,\end{equation}
і що функція $h:\overline{{\Bbb R}^n}\times \overline{{\Bbb
R}^n}\rightarrow {\Bbb R}$ задовольняє умову слабкої сферикалізації.

\medskip
Справедливе наступне твердження.

\medskip
\begin{theorem}\label{th1}
{\sl\, Нехай $n\geqslant 2,$ $D$ -- область в ${\Bbb M}^n,$ $x_0\in
D,$ $a, b\in {\Bbb M}^n_*,$ $a\ne b$ і нехай
$f:D\setminus\{x_0\}\rightarrow {\Bbb M}_*^n\setminus\{a, b\}$ --
відкрите дискретне відображення, яке задовольняє
умови~(\ref{eq2*A})--(\ref{eqA2}) у всіх точках $y_0\in
\overline{D^{\,\prime}},$ де $D^{\,\prime}:=f(D\setminus\{x_0\}).$
Припустимо, що многовид ${\Bbb M}^n_*$ є зв'язним і допускає умову
слабкої сферикалізації. Припустимо, крім того, що для кожної точки
$y_0\in \overline{D^{\,\prime}}$ знайдеться $r_0=r_0(y_0)>0$ таке,
що $q_{y_0}(r)<\infty$ при майже всіх $r\in (0, r_0).$ Тоді $f$ має
неперервне продовження $\overline{f}:D\rightarrow\overline{{\Bbb
M}^n},$ неперервність якого слід розуміти в сенсі метрики $h.$
Продовжене відображення є відкритим і дискретним у $D.$

Зокрема, твердження теореми~\ref{th1} виконується, якщо $Q\in
L^1(D^{\,\prime}).$ В цьому випадку, якщо $f(x_0)\ne \infty,$ то
знайдеться окіл $U=U(x_0, f)\subset D$ точки $x_0,$ залежний тільки
від $x_0$ та відображення $f,$ і стала $C=C(n, D, x_0)>0$ такі, що
\begin{equation}\label{eq2C}
d_*(\overline{f}(x), \overline{f}(x_0))\leqslant\frac{C_n\cdot
(\Vert Q\Vert_1)^{1/n}}{\log^{1/n}\left(1+\frac{\delta}{2d(x,
x_0)}\right)}
\end{equation}
для всіх $x, y\in U,$
де $\Vert Q\Vert_1$ -- норма функції $Q$ в $L^1(D^{\,\prime}).$ }
\end{theorem}

\medskip
З теореми~\ref{th1} миттєво отримаємо наступний

\medskip
\begin{corollary}
{\sl\, Нехай $n\geqslant 2,$ $D$ -- область в ${\Bbb R}^n,$ $x_0\in
D,$ $a, b\in {\Bbb R}^n,$ $a\ne b$ і нехай
$f:D\setminus\{x_0\}\rightarrow {\Bbb R}^n\setminus\{a, b\}$ --
відкрите дискретне відображення, яке задовольняє
умови~(\ref{eq2*A})--(\ref{eqA2}) у всіх точках $y_0\in
\overline{D^{\,\prime}},$ де $D^{\,\prime}:=f(D\setminus\{x_0\}).$
Припустимо, що для кожної точки $y_0\in \overline{D^{\,\prime}}$
знайдеться $r_0=r_0(y_0)>0$ таке, що $q_{y_0}(r)<\infty$ при майже
всіх $r\in (0, r_0).$ Тоді $f$ має неперервне продовження
$\overline{f}:D\rightarrow\overline{{\Bbb R}^n},$ неперервність
якого слід розуміти в сенсі хордальної метрики $h.$ Продовжене
відображення є відкритим і дискретним у $D.$

Зокрема, твердження теореми~\ref{th1} виконується, якщо $Q\in
L^1(D^{\,\prime}).$ В цьому випадку, якщо $f(x_0)\ne \infty,$ то
знайдеться окіл $U=U(x_0, f)\subset D$ точки $x_0,$ залежний тільки
від $x_0$ та відображення $f,$ і стала $C=C(n, D, x_0)>0$ такі, що
\begin{equation}\label{eq2G}
|\overline{f}(x)-\overline{f}(x_0)|\leqslant\frac{C_n\cdot (\Vert
Q\Vert_1)^{1/n}}{\log^{1/n}\left(1+\frac{\delta}{2|x-x_0|}\right)}
\end{equation}
для всіх $x, y\in U,$
де $\Vert Q\Vert_1$ -- норма функції $Q$ в $L^1(D^{\,\prime}).$  }
\end{corollary}

\medskip
\begin{remark}\label{rem2}
Відзначимо, що усі квазірегулярні відображення
$f:D\setminus\{x_0\}\rightarrow {\Bbb R}^n$ задовольнять умову
\begin{equation}\label{eq22}
M(\Gamma_f(y_0, r_1, r_2))\leqslant \int\limits_{f(D)\cap
A(y_0,r_1,r_2)} K_O\cdot N(y, f, D\setminus\{x_0\})\cdot \eta^n
(|y-y_0|)\, dm(y) \end{equation}
у кожній точці $y_0\in \overline{f(D)}\setminus\{\infty\}$ з деякою
сталою $K_O=K_O(f)\geqslant 1$ і довільною вимірною за Лебегом
функцією $\eta: (r_1,r_2)\rightarrow [0,\infty],$ яка задовольняє
умову~(\ref{eqA2}). Тут функція кратності $N(y, f,
D\setminus\{x_0\})$ визначена у~(\ref{eq23}).

Справді, квазірегулярні відображення задовольняють умову
\begin{equation}\label{eq24}
M(\Gamma_f(y_0, r_1, r_2))\leqslant \int\limits_{f(D)\cap
A(y_0,r_1,r_2)} K_O\cdot N(y, f, D\setminus\{x_0\})\cdot
(\rho^{\,\prime})^n(y)\, dm(y) \end{equation}
для довільної функції $\rho^{\,\prime}\in{\rm adm}\, f(\Gamma_f(y_0,
r_1, r_2)),$ див. \cite[зауваження~2.5.II]{Ri}. Покладемо тепер
$\rho^{\,\prime}(y):=\eta(|y-y_0|)$ для $y\in A(y_0,r_1,r_2)\cap
f(D),$ і $\rho^{\,\prime}(y)=0$ у решті точок. З огляду на теорему
Лузіна, можна припустити, що функція $\rho^{\,\prime}$ є вимірною за
Борелем (див., напр., \cite[розділ~2.3.6]{Fe}). Тоді, з огляду
на~\cite[теорема~5.7]{Va},
$$\int\limits_{\gamma_*}\rho^{\,\prime}(y)\,|dy|\geqslant
\int\limits_{r_1}^{r_2}\eta(r)\,dr\geqslant 1$$
для кожної (локально спрямлюваної) кривої $\gamma_*\in \Gamma(S(y_0,
r_1), S(y_0, r_2), A(y_0, r_1, r_2)).$ Підставляючи обрану функцію
$\rho^{\,\prime}$ у~(\ref{eq24}), ми отримаємо бажане
співвідношення~(\ref{eq22}).
\end{remark}

\medskip
Враховуючи зауваження~\ref{rem2}, з огляду на теорему~\ref{th1},
маємо наступний

\medskip
\begin{corollary}\label{cor1}
{\sl\, Нехай $n\geqslant 2,$ $D$ -- область в ${\Bbb R}^n,$ $x_0\in
D,$ $f:D\setminus\{x_0\}\rightarrow {\Bbb R}^n\setminus\{a, b\},$
$a, b\in {\Bbb R}^n,$ $a\ne b$ -- квазірегулярне відображення. Якщо
$N(y, f, D\setminus\{x_0\})\in L^1(D^{\,\prime}),$ то $f$ має
неперервне продовження $\overline{f}:D\rightarrow\overline{{\Bbb
R}^n},$ неперервність якого слід розуміти в сенсі хордальної метрики
$h$ у~(\ref{eqA17}).}
\end{corollary}

\medskip
На жаль, твердження наслідку~\ref{cor1} не є новим. Воно випливає з
теореми Сохоцького-Вейєрштрасса, яка стверджує, що для істотно
особливих точок $x_0$ відображення $f$ виконана умова: $N(y, f,
D\setminus\{x_0\})=\infty$ для всіх $y\in\overline{{\Bbb
R}^n}\setminus E,$ де $E$ -- множина ємності нуль (див.
\cite[наслідок~2.11.III]{Ri}). Отже, для істотно особливих точок
$x_0$ відображення $f$ умова $N(y, f, D\setminus\{x_0\})\in
L^1(D^{\,\prime}),$ присутня у наслідку~\ref{cor1}, a priori не може
виконуватися.

\medskip
Тим не менш, само доведення теореми~\ref{th1} (а отже, і
наслідку~\ref{cor1}) відбувається за принципово новою схемою, яку
наразі не можна вважати загальновідомою в теорії відображень.

\medskip
{\bf 2. Допоміжні твердження.} Перед початком доведення основних
результатів встановимо наступні твердження.

\medskip
\begin{proposition}\label{pr1}
{\sl\, Нехай $D$ -- область в ${\Bbb M}^n,$ $n\geqslant 2,$ $x_0\in
D,$ $f:D\setminus\{x_0\}\rightarrow {\Bbb M}_*^n,$ $n\geqslant 2,$
-- довільне відображення. Припустимо, що многовид ${\Bbb M}^n_*$ є
зв'язним і допускає умову слабкої сферікалізації. Тоді $C(x_0, f)$
-- непорожній континуум в $\overline{{\Bbb M}_*^n}.$ }
\end{proposition}

Тут і надалі
$$C(x, f):=\{y\in \overline{{\Bbb R}^n}:\exists\,x_k\in D:
x_k\rightarrow x, f(x_k) \rightarrow y, k\rightarrow\infty\}\,.$$
\begin{proof}
Покладемо $r_0=d(z_0, \partial D).$ Нехай $r_m>0,$ $m\in {\Bbb N}$
-- довільна послідовність така, що $r_m\rightarrow 0$ при
$m\rightarrow \infty.$ Нехай $m_0\in {\Bbb N}$ таке, що $r_m<r_0$
при всіх $m>m_0.$ Зауважимо, що
$C(x_0, f)=\bigcap\limits_{m=1}^{\infty}\overline{f(B(x_0,
r_m)\setminus\{x_0\})}.$
Звідси випливає, що $C(x_0, f)$ замкнена множина як зчисленний
перетин замкнених множин. Крім того, $C(x_0, f)$ -- континуум як
перетин спадної послідовності континуумів (див.
\cite[теорема~5.II.47.5]{Ku}). Оскільки простір $\overline{{\Bbb
M}_*^n}$ компактний, $C(x_0, f)$ не є порожнім з огляду на умову
Кантора щодо непорожнього перетину довільної системи замкнених
спадаючих множин (див. \cite[$\S\,41,$ гл.~4, пункт~I]{Ku}).~$\Box$
\end{proof}

\medskip
Сформулюємо дуже важливе топологічне твердження, яке неодноразово
буде використовуватися в подальшому (див.
\cite[те\-о\-ре\-ма~1.I.5.46]{Ku}).

\medskip
\begin{proposition}\label{pr2}
{\sl\, Нехай $A$ -- множина в топологічному просторі $X.$ Якщо
множина $C$ є зв'язною і $C\cap A\ne \varnothing\ne C\setminus A,$
то $C\cap
\partial A\ne\varnothing.$}
\end{proposition}

\medskip
Нехай $D\subset {\Bbb M}^n,$ $f:D\rightarrow {\Bbb M}_*^n$ є
дискретним відкритим відображенням, $\beta: [a,\,b)\rightarrow {\Bbb
R}^n$ -- крива, і $x\in\,f^{\,-1}(\beta(a)).$ Крива $\alpha:
[a,\,c)\rightarrow D$ називається {\it максимальним підняттям кривої
$\beta$ при відображенні $f$ з початком у точці $x,$} якщо $(1)\quad
\alpha(a)=x\,;$ $(2)\quad f\circ\alpha=\beta|_{[a,\,c)};$ $(3)$\quad
яким би не було $c<c^{\,\prime}\leqslant b,$ не існує кривої
$\alpha^{\prime}: [a,\,c^{\prime})\rightarrow D$ такої, що
$\alpha=\alpha^{\prime}|_{[a,\,c)}$ і $f\circ
\alpha^{\,\prime}=\beta|_{[a,\,c^{\prime})}.$ Аналогічно можна
визначити максимальне підняття $\alpha: (c,\,b]\rightarrow D$ кривої
$\beta: (a,\,b]\rightarrow {\Bbb M}_*^n$ при відображені $f$ з
кінцем у точці $x\in\,f^{\,-1}(\beta(b)).$ Максимальне підняття
$\alpha: [a,\,c)\rightarrow D$ кривої $\beta: [a,\,b)\rightarrow
{\Bbb M}_*^n$ при відображенні $f$ з початком у точці $x$
називається {\it повним,} якщо в наведеному вище означенні $c=b.$

\medskip
Якщо $\beta:[a, b)\rightarrow\overline{{\Bbb M}_*^n}$ -- крива і
$C\subset {\Bbb M}^n,$ то ми будемо говорити, що $\beta\rightarrow
C$ при $t\rightarrow b,$ якщо $d(\beta(t), C)\rightarrow 0$ при
$t\rightarrow b$ (див. \cite[розділ~3.11]{MRV$_3$}), де $d(\beta(t),
C)$ визначено в~(\ref{eq1A}). Справедливо наступне твердження
(див.~\cite[лема~3.12]{MRV$_3$} для випадку ${\Bbb M}^n={\Bbb
M}_*^n={\Bbb R}^n$).

\medskip
\begin{proposition}\label{pr3}
{\sl Нехай многовид ${\Bbb M}^n$ є зв'язним, $f:D\rightarrow {\Bbb
M}^n,$ $n\geqslant 2,$ $D\ne {\Bbb M}^n,$ -- відкрите дискретне
відображення, нехай $x_0\in D,$ $\overline{D}$ -- компакт у ${\Bbb
M}^n,$ і нехай $\beta: [a,\,b)\rightarrow {\Bbb M}_*^n$ --  крива,
така що $\beta(a)=f(x_0)$ і, крім того, або існує
$\lim\limits_{t\rightarrow b}\beta(t),$ або $\beta(t)\rightarrow
\partial f(D)$ при $t\rightarrow b.$ Тоді $\beta$ має максимальне
підняття $\alpha: [a,\,c)\rightarrow D$ при відображенні $f$ з
початком у точці $x_0.$ Більше того, якщо $\beta(t)\rightarrow
\partial D^{\,\prime}$ при $t\rightarrow b-0,$ то будь-яке її максимальне підняття
$\alpha:[a,c)\rightarrow D$ з початком в довільній точці $x\in
f^{\,-1}(\beta(a))$ при відображенні $f$ задовольняє умову:
$d(\alpha(t),\partial D)\rightarrow 0$ при $t\rightarrow c-0.$}
\end{proposition}

\medskip
\begin{proof}
Існування максимальних піднять $\alpha: [a,\,c)\rightarrow D$ при
відображенні $f$ з початком у точці $x_0$ доведено в~\cite{SM}.
Доведемо другу частину твердження, для чого будемо користуватися
методологією доведення~\cite[лема~3.12]{MRV$_3$}, див. також
\cite[лема~2.1]{Skv}.

Передусім зауважимо, що формулювання твердження~\ref{pr3} є
коректним, а саме, за умов твердження $\partial D$ не порожня.
Дійсно, оскільки за умовою многовид ${\Bbb M}^n$ є зв'язним і $D\ne
{\Bbb M}^n,$ можна з'єднати точки $x_1\in D$ і $x_2\in D\setminus
{\Bbb M}^n$ кривою $\gamma:[0, 1]\rightarrow {\Bbb M}^n,$
$\gamma(0)=x_1$ і $\gamma(1)=x_2.$ Тоді $|\gamma|\cap
D\ne\varnothing\ne {\Bbb M}^n\setminus D,$ тому за твердженням
\ref{pr2} $\gamma\cap
\partial D\ne\varnothing.$ Отже, $\partial D\ne\varnothing,$ що і треба було довести.

\medskip
Оскільки $\partial D\ne\varnothing,$ то можливі два випадки: або
$d(\alpha(t),\partial D)\rightarrow 0$ при $t\rightarrow c-0,$ або
існує $\delta_0>0$ і послідовність $t_k\rightarrow c-0$ при
$k\rightarrow\infty$ така, що $d(\alpha(t_k),\partial
D)\geqslant\delta_0>0$ при $k\rightarrow\infty.$

\medskip
Доведемо від супротивного, що друга можливість реалізуватися не
може. Припустимо протилежне. Тоді розглянемо множину
$$D_0=\left\{x\in {\Bbb M}^n:  x=\lim_{k\rightarrow \infty} \alpha(t_k),\quad
t_k \in [a,c), \quad\lim\limits_{k\rightarrow \infty}
t_k=c\right\}.$$
Доведемо спочатку, що $c\ne b.$ Дійсно, послідовність точок
$x_k:=\alpha(t_k)$ можна вважати збіжною до деякої точки $x_0\in D$
з огляду на те, що $\overline{D}$ -- компакт в ${\Bbb M}^n.$ Тоді,
оскільки $f$ є неперервним у $D,$ крім того, $f$ є відкритим
відображенням, то $f(x_k)=\beta(t_k)\rightarrow y_0\in D^{\,\prime}$
при $k\rightarrow\infty,$ проте це суперечить означенню кривої
$\beta.$

\medskip
Отже, $c\ne b.$ Нехай $x\in D_0\cap D,$ тоді з огляду на
неперервність $f$ ми отримаємо, що $f(\alpha(t_k))\rightarrow f(x)$
при $k\rightarrow \infty,$ де $t_k \in [a,c), t_k\rightarrow c$ при
$k\rightarrow \infty.$ Однак, $f(\alpha(t_k))=\beta(t_k)\rightarrow
\beta(c)$ при $k\rightarrow \infty.$ Тоді відображення $f$ стале на
$D_0\cap D.$ З іншого боку, множина $\overline{|\alpha|}$ є
компактною, оскільки $\overline{|\alpha|}$ -- замкнена підмножина
компактного простору $\overline{D}$ (див.~\cite[теорема~2.II.4, \S
41]{Ku}). За умовою Кантора на компакті $\overline{|\alpha|}$ та з
огляду на монотонність послідовності зв'язних множин
$\alpha([t_k,c)),$
$D_0=\bigcap\limits_{k=1}^{\infty} \overline{\alpha([t_k,c))}
\ne\varnothing$
(див.~\cite[пункт~(2').I.41.4]{Ku}). З огляду
на~\cite[теорема~5.II.47.5]{Ku} множина $D_0$ зв'язна.

В такому випадку, $D_0$ складається з однієї точки. Дійсно, якщо
$D_0$ -- невироджений континуум то, оскільки за побудовою $x_0\in
D_0\cap D,$ то маємо й невироджений підконтинуум $E_0\subset D\cap
D_0.$ Останнє суперечить дискретності відображення $f$ у $D,$ бо
тоді з огляду на неперервність відображення $f$ ми мали б, що
$f(E_0)=\{\beta(c)\}.$ Отже, $D_0=\{x_0\}.$

\medskip
У такому випадку, криву $\alpha:[a,c) \rightarrow D$ можна
продовжити до замкненої кривої $\alpha:[a,c] \rightarrow D,$ причому
$f(\alpha(c))=\beta(c).$ Тоді за~\cite[лема~2.1]{SM} знайдеться ще
одне максимальне підняття $\alpha^{\prime}$ кривої $\beta|_{[c,b)}$
з початком в точці $\alpha(c).$ Об'єднуючи підняття $\alpha$ i
$\alpha^{\prime},$ ми отримаємо нове підняття
$\alpha^{\prime\prime}$ кривої $\beta,$ визначене на $[a,
c^{\prime}),$ $c^{\,\prime}\in (c,b),$ що суперечить
<<максимальності>>\, вихідного підняття $\alpha.$ Отримане
протиріччя вказує на те, що припущення $d(\alpha(t_k),\partial
D)\geqslant\delta_0>0$ при деяких послідовності $t_k\rightarrow
c-0,$ $k\rightarrow\infty,$ і $\delta_0>0$ було невірним.
Твердження~\ref{pr3} повністю доведено.
\end{proof}

\medskip
Є справедливим наступний результат.

\medskip
\begin{proposition}\label{pr4}
{\sl\,Нехай $n\geqslant 2,$ $D$ -- область в ${\Bbb R}^n,$ $x_0\in
D,$ $f:D\setminus\{x_0\}\rightarrow {\Bbb R}^n$ -- відкрите
дискретне відображення. Якщо $f$ має неперервне продовження
$\overline{f}:D\rightarrow \overline{{\Bbb R}^n}$ у точку $x_0,$ то
$\overline{f}$ також є відкритим і дискретним відображенням. }
\end{proposition}

\medskip
\begin{proof}
Без обмеження загальності, застосовуючи за потреби додаткове
перетворення інверсії $\psi(x)=\frac{x}{|x|^2},$ $\infty\mapsto 0,$
ми можемо вважати, що $\overline{f}(x_0)\ne\infty.$

Усюди далі $\mu\,(y,\,f,\,G)$ -- топологічний індекс відображення
$f$ у точці  $y\,\in\,f(G)\setminus f(\partial G)$ відносно області
$G\subset D,$ $\overline{G}\subset D.$ Будемо говорити, що
ві\-доб\-ра\-жен\-ня $f$ {\it зберігає орієнтацію,} якщо
топологічний індекс $\mu\, (y,\,f,\,G)$ задовольняє умову $\mu\,
(y,\,f,\,G)\,>\,0$ для довільної області $G\subset D,$
$\overline{G}\subset D,$ і $y\,\in\,f(G)\setminus f(\partial G),$
див., напр., п.~4 розд.~I, с.~17 \cite{Ri}.

Відомо, що дискретні відкриті ві\-доб\-ра\-жен\-ня в ${\Bbb R}^n,$
$n\geqslant 2,$ або зберігають орієнтацію, або анти-зберігають,
див., напр., п. 4 розд. I в \cite{Ri}. Нехай $f,$ наприклад,
зберігає орієнтацію. Покажемо, що продовжене ві\-доб\-ра\-жен\-ня
$\overline{f}$ зберігає орієнтацію, відкрито і дискретно. Позначимо,
як звично, через $B_f\left(D\setminus\{x_0\}\right)$ множину точок
розгалуження ві\-доб\-ра\-жен\-ня $f$ в області $D\setminus\{x_0\},$
а через $B_{\overline{f}}(D)$ -- множину точок розгалуження
продовженого у точку $x_0$ ві\-доб\-ра\-жен\-ня $\overline{f}$ в
області $D.$ Якщо $x_0$ -- точка локальної гомеоморфності
ві\-доб\-ра\-жен\-ня $\overline{f},$ доводити нема що. Нехай точка
$x_0\,\in\, B_{\overline{f}}(D).$ За теоремою Чернавського ${\rm
dim}\,B_f(D\setminus\{x_0\})=\,{\rm
dim}\,\overline{f}(B_f(D\setminus\{x_0\}))\,\leqslant\,n-2,$ див.,
напр., теорему~4.6 розд.~I \cite{Ri}, де ${\rm dim}$ позначає
топологічну розмірність множини, див. \cite{HW}. Тоді отримаємо:
\begin{equation}\label{eq38*!}
{\rm dim}\,f(B_{\overline{f}}(D))\leqslant n-2\,,
\end{equation}
оскільки
$f(B_{\overline{f}}(D))\,=\,f(B_{\overline{f}}(D\setminus\{x_0\}))\bigcup
\left\{\overline{f}(x_0)\right\},$ множина
$\left\{\overline{f}(x_0)\right\}$ є замкненою і топологічна
розмірність кожної з множин $\overline{f}(B_f(D\setminus\{x_0\}))$ і
$\left\{\overline{f}(x_0)\right\}$ не перевищує $n-2,$ див.
наслідок~1 розд.~III п.~3 У \cite{HW}. Нехай $G$ -- область в $D$
така, що $\overline{G}\subset D$ і $y\in\, f(G)\setminus f(\partial
G).$ Тоді, в силу (\ref{eq38*!}), існує точка
$y_0\,\notin\,f(B_{\overline{f}}(D)),$ що належить до тієї ж
компоненти зв'язності множини $\overline{{\Bbb R}^n}\setminus
f(\partial G),$ що і $y.$

\medskip
Нагадаємо, що для довільного ві\-доб\-ра\-жен\-ня $f:D\rightarrow
{\Bbb R}^n$ існує область $G\subset \overline{D}$ така, що
$\overline{G}\cap f^{\,-1}(f(x))=\{x\}.$ Тоді величина
$\mu\,(f(x),\,f,\,G),$ яка називається {\it локальним топологічним
індексом}, не залежить від обирання області $G$ і позначається
символом $i(x,f).$

\medskip
В силу того, що топологічний індекс $\mu(y, f , G)$ є сталою
величиною на кожній компоненті зв'язності множини $\overline{{{\Bbb
R}^n}}\setminus f(\partial G)$ (див. властивість~$D_1$
у~\cite[пропозиція~4.4.I]{Ri}), будемо мати:
$\mu(y, \overline{f}, G)\,=\,\mu(y_0, f, G)\,=
\,\sum\limits_{x\,\in\,G\cap \{f^{\,-1}(y_0)\}} \,i(x,f)\,>\,0.$
Отже, ві\-доб\-ра\-жен\-ня $f$ зберігає орієнтацію в $D.$ Нарешті,
щодо будь--якого $y\,\in f(D),$ в силу дискретності
ві\-доб\-ра\-жен\-ня $f$ в області $D\setminus\{x_0\},$ множина
$\left\{f^{\,-1}(y)\right\}$ є не більш ніж зчисленною і тому ${\rm
dim}\,\left\{f^{\,-1}(y)\right\}=0.$ Отже, за~\cite{TY}, с.~333,
ві\-доб\-ра\-жен\-ня $f$ є відкритим і дискретним, що і потрібно
було довести.~$\Box$
\end{proof}

\medskip
Наступне твердження містить в собі елементарний зв'язок між
інтегровністю функції по сферах та інтегровністю по мірі Лебега, і є
однією з версій теореми Фубіні (див. співвідношення~(3.3)--(3.4) у
\cite{IS}).

\medskip
\begin{proposition}\label{pr5}
{\sl\, Нехай $D$ -- область в ${\Bbb M}^n,$ $n\geqslant 2,$
$Q:D^{\,\prime}\rightarrow [0, \infty],$ $Q\equiv 0$ на ${\Bbb
R}^n\setminus D^{\,\prime}.$ Тоді для всякого $y_0\in
\overline{D^{\,\prime}}\setminus\{\infty\}$ знайдуться нормальний
окіл точки $y_0$ і сталі $C_1=C_1(y_0), C_2=C_2(y_0)>0$ такі, що для
всяких
$0\leqslant\varepsilon_1<\varepsilon_2<\operatorname{dist}(x_0,\partial
U)$ виконуються нерівності
\begin{equation}\label{eq1U}
C_1\int\limits_{\varepsilon_1}^{\varepsilon_2} \int\limits_{S(y_0,
r)}Q(y)\,d\mathcal{A}(y)dr\leqslant
\int\limits_{\varepsilon_1<|y-y_0|<\varepsilon_2}Q(y)\,dm(y)\leqslant
C_2\int\limits_{\varepsilon_1}^{\varepsilon_2} \int\limits_{S(y_0,
r)}Q(y)\,d\mathcal{A}(y)dr\,.
\end{equation}
Співвідношення~(\ref{eq1U}) передбачає існування (скінченного або
нескінченного) інтегралу від функції $Q$ по сфері $S(y_0, r)$ при
майже кожному $\varepsilon_1\leqslant r\leqslant \varepsilon_2$
відносно міри $\mathcal{A}$ поверхні $S(y_0, r).$}
\end{proposition}

\medskip
Для області $D\subset {\Bbb R}^n,$ $n\geqslant 2,$ і вимірної за
Лебегом функції $Q:{\Bbb R}^n\rightarrow [0, \infty]$ визначимо
через $\frak{F}_Q(D)$ сім'ю всіх відкритих дискретних відображень
$f:D\rightarrow {\Bbb R}^n$ таких, що співвідношення~(\ref{eq2*A})
виконується при $p=n$ для кожної точки $y_0\in f(D).$ Є справедливим
наступне твердження (див., напр., \cite[теорема~1.1]{SSD}).

\medskip
\begin{proposition}\label{pr6}
{\sl\, Нехай $Q\in L^1({\Bbb R}^n).$ Тоді знайдеться стала $C_n>0,$
яка залежить тільки від розмірності простору $n,$ така що для
будь-якого $x_0\in D$ і будь-якого $r_0>0$ такого, що $0<2r_0<{\rm
dist}\,(x_0,
\partial D),$ виконується нерівність
\begin{equation}\label{eq2CB}
|f(x)-f(x_0)|\leqslant\frac{C_n\cdot (\Vert
Q\Vert_1)^{1/n}}{\log^{1/n}\left(1+\frac{r_0}{|x-x_0|}\right)}
\quad\forall\,\,x\in B(x_0, r_0)\,,\quad \forall\,\, f\in
\frak{F}_Q(D)\,,
\end{equation}
де $\Vert Q\Vert_1$ -- норма функції $Q$ в $L^1({\Bbb R}^n).$
Зокрема, сім'я $\frak{F}_Q(D)$ є одностайно неперервною в $D.$}
\end{proposition}

\medskip
Має місце наступне твердження, див. \cite[лема~2.2]{IS$_2$}, див.
також~\cite[теорема~10.12]{Va} і \cite[лема~2.2]{SevSkv$_2$}.

\medskip
\begin{proposition}\label{pr7}
{\sl\, Нехай $D$ -- область в ${\Bbb M}^n,$ $n\geqslant 2,$ і
$x_0\in D.$ Тоді для кожного $P>0$ і будь-якого околу $U$ точки
$x_0$ знайдеться окіл $V\subset U$ цієї ж точки, такий що нерівність
$M(\Gamma(E, F, D))>P$ виконується для довільних континуумів $E,
F\subset D,$ що перетинають $\partial U$ и $\partial V.$}
\end{proposition}

\medskip
\begin{remark}\label{rem1}
Оскільки модуль сім'ї кривих, що проходять через фіксовану точку,
дорівнює нулю (див. пункт~7.9 у \cite{Va}), то твердження~\ref{pr7}
залишається вірним для випадку, коли $x_0$ -- ізольована точка межі
області $D.$
\end{remark}

\medskip
Наступне твердження доведено в~\cite[лема~2.1]{IS$_2$}.

\medskip
\begin{proposition}\label{pr8}
{\sl Нехай $a,$ $b,$ $c$ і $d$
--- чотири різні точки області $D$ ріманового многовиду
${\Bbb M}^n,$ $n\geqslant 2.$ Тоді знайдуться непересічні жорданові
криві $\gamma_1\colon[0,1]\rightarrow D$ і $\gamma_2\colon[0,
1]\rightarrow D,$ такі, що $\gamma_1(0)=a,$ $\gamma_1(1)=b,$
$\gamma_2(0)=c$ і $\gamma_2(1)=d.$}
\end{proposition}

\medskip
{\bf 3. Доведення теореми~\ref{th1}.} Доведення розіб'ємо на
наступні кроки.

\medskip
{\bf I. Існування неперервного продовження в точку $x_0.$} Частково
скористаємося і удосконалимо підхід, здійснений при доведенні
теореми~1 в \cite{Sev$_3$}, див. також теорему~6 в \cite{Sev$_2$}.
Припустимо супротивне, а саме, що відображення $f$ не має границі в
точці $x_0$ в сенсі ``хордальної'' метрики $h.$ З огляду на
компактність простору $\overline{{\Bbb M}^n}$ знайдуться $z_1,
z_2\in C(x_0, f),$ $z_1\ne z_2.$ За твердженням~\ref{pr1} $C(x_0,
f)$ є невиродженим континуумом, і оскільки простір $\overline{{\Bbb
M}^n}$ є метричним, цей континуум містить безліч точок (див., напр.,
\cite[теорема~4.I.5, \S\,46]{Ku}). Отже, можна вважати, що
$z_1\ne \infty\ne z_2.$
Тоді також знайдуться послідовності $x_m, x_m^{\,\prime}\in
D\setminus \{x_0\},$ $m=1,2,\ldots,$ такі що $x_m,
x_m^{\,\prime}\rightarrow x_0$ при $m\rightarrow\infty,$ причому
$y_m:=f(x_m)\rightarrow z_1$ і
$y^{\,\prime}_m:=f(x_m^{\,\prime})\rightarrow z_2$ при
$m\rightarrow\infty.$ Оскільки за означенням слабкої сферикалізації
метрики $d$ і $h$ формують одну і ту саму топологію на ${\Bbb M}^n,$
вказану збіжність можна розуміти в сенсі метрики $d.$

\begin{figure}
\centerline{\includegraphics[scale=0.6]{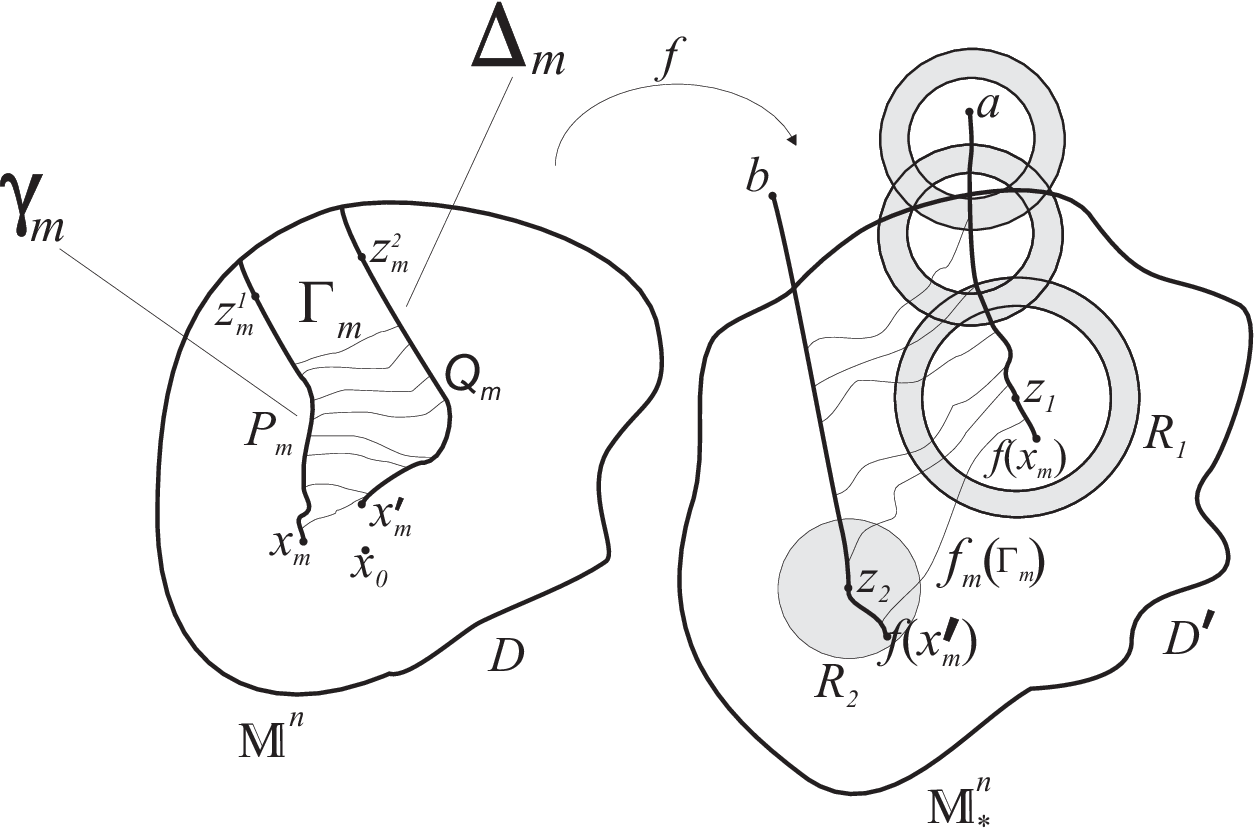}} \caption{До
доведення теореми~\ref{th1}. Випадок~1)}\label{figure1}
 \end{figure}

\medskip
За твердженням~\ref{pr8} точки $a$ та $z_1$ та $b$ та $z_2$ можна
з'єднати кривими $\gamma_1\colon [1/2, 1] \rightarrow {\Bbb M}^n_*$
і $\gamma_2\colon[1/2, 1] \rightarrow {\Bbb M}^n_*.$ Якщо $a=z_1,$
або $b=z_2,$ відповідні криві $\gamma_1$ та/або $\gamma_2$
вважаються виродженими (одноточковими). Оберемо $R_1>0$ та $R_2>0$
такими, що
$$
\left(\overline{B(z_1, R_1)}\cup
|\gamma_1|\right)\cap\left(\overline{B(z_2, R_2)}\cup
|\gamma_2|\right)=\varnothing\,.
 $$
Можна вважати, що $y_m:=f(x_m)\in B(z_1, R_1)$ і
$y^{\,\prime}_m:=f(x_m^{\,\prime})\in B(z_2, R_2)$ при всіх $m\in
{\Bbb N}.$ З'єднаємо точки $f(x_m)$ і $z_1$ кривою
$\alpha^{\,*}_m\colon [0, 1/2] \rightarrow B(z_1, R_1),$ і
$f(x^{\,\prime}_m)$ і $z_2$ кривою $\beta^{\,*}_m\colon[0,
1/2]\rightarrow B(z_2, R_2)$ (див. малюнок~\ref{figure1}). Це є
можливим, оскільки інфінітезимальні кулі на многовиді є зв'язними.
Покладемо
$$
\alpha_m(t)=\quad\left\{\begin{array}{rr}
\alpha^*_m(t), & t\in [0, 1/2],\\
\gamma_1(t), & t\in [1/2, 1]\end{array} \right.\,,\quad
\beta_m(t)=\quad\left\{
\begin{array}{rr}
\beta^*_m(t), & t\in [0, 1/2],\\
\gamma_2(t), & t\in [1/2, 1].\end{array} \right.
 $$
За побудовою множини
$A_1:=|\gamma_1|\cup \overline{B(z_1, R_1)},$ $A_2:=|\gamma_2|\cup
\overline{B(z_2, R_2)}$
не перетинаються, зокрема, знайдеться $\varepsilon^{\,*}_1>0$ таке,
що
\begin{equation}\label{eq1J}
d(A_1,A_2)\geqslant \varepsilon^{\,*}_1>0\,.
\end{equation}
З огляду на визначення кривих $\alpha_m$ та $\beta_m$
$$|\alpha_m|\cap D^{\,\prime}\ne\varnothing \ne
|\alpha_m|\cap\left({\Bbb M}_*^n\setminus D^{\,\prime}\right)\,,$$
$$|\beta_m|\cap D^{\,\prime}\ne\varnothing \ne
|\beta_m|\cap\left({\Bbb M}_*^n\setminus D^{\,\prime}\right)\,.$$
Тоді з огляду на твердження~\ref{pr2} існують точки $t_1,
t_2\in\left(\frac{1}{2}, 1\right)$ такі, що
$\alpha_m(t_1)\in\partial D^{\,\prime},$ $\beta_m(t_2)\in\partial
D^{\,\prime}.$
Без обмеження загальності можна вважати, що $\alpha_m(t)\in
D^{\,\prime}$ при всіх $t\in [0, t_1)$ і $\beta_m(t)\in
D^{\,\prime}$ при всіх $t\in [0, t_2).$ Покладемо
$$\widetilde{\alpha}_m:=\alpha_m|_{[0, t_1)}\,,
\quad \widetilde{\beta}_m:=\beta_m|_{[0, t_2)}\,.$$
Розглянемо максимальні підняття $\gamma_m:[0, c_m)\rightarrow
D\setminus\{x_0\}$ і $\Delta_m:[0, d_m)\rightarrow
D\setminus\{x_0\}$ кривих $\widetilde{\alpha}_m$ і
$\widetilde{\beta}_m$ в області $D\setminus\{x_0\}$ з початками у
точках $x_m$ та $x^{\,\prime}_m.$ Такі підняття існують з огляду на
твердження~\ref{pr3}. Завдяки цьому ж твердженню
$\gamma_m(t)\rightarrow \partial (D\setminus\{x_0\})$ при
$t\rightarrow c_m-0$ і $\Delta_m(t)\rightarrow \partial
(D\setminus\{x_0\})$ при $t\rightarrow d_m-0.$ Можливі 4 наступні
ситуації:

\medskip
\textbf{1)} при всіх $m\in {\Bbb N},$ $\gamma_m(t)\rightarrow
\partial D$ при $t\rightarrow c_m-0$ і $\Delta_m(t)\rightarrow
\partial D$ при $t\rightarrow d_m-0;$

\medskip
\textbf{2)} при всіх $m\in {\Bbb N},$ $\gamma_m(t)\rightarrow
\partial D$ при $t\rightarrow c_m-0,$ але існує $m_0\in {\Bbb N}$
такий, що $\Delta_{m_0}(t)\rightarrow 0$ при $t\rightarrow
d_{m_0}-0;$

\medskip
\textbf{3)} існує $k_0\in {\Bbb N}$ такий, що
$\gamma_{k_0}(t)\rightarrow 0$ при $t\rightarrow c_{k_0}-0,$ проте,
при всіх $m\in {\Bbb N},$ $\Delta_m(t)\rightarrow \partial D$ при
$t\rightarrow d_m-0;$

\medskip
\textbf{4)} існують $k_0, m_0\in {\Bbb N}$ такі, що
$\gamma_{k_0}(t)\rightarrow 0$ при $t\rightarrow c_{k_0}-0$ і
$\Delta_{m_0}(t)\rightarrow 0$ при $t\rightarrow d_{m_0}-0.$

\medskip
У подальшому нам слід розглянути кожну з цих ситуацій в окремості,
причому з двох ситуацій 2), або 3), достатньо розглянути тільки одну
(бо ці ситуації, як легко бачити, відрізняються заміною позначень).

\medskip
Розглянемо \textbf{випадок~1).} Нехай $r_0=r_0(y)>0$ -- число з умов
теореми, яке визначено для кожного $y_0\in\overline{D^{\,\prime}}.$
За твердженням~\ref{pr5} для кожного $y_0\in{\Bbb M}^n_*$ знайдеться
$\delta(y_0)>0$ і стала $C=C(y_0)>0$ такі, що
\begin{equation}\label{eq7A}
\int\limits_{\varepsilon_1<d_*(y,
y_0)<\varepsilon_2}Q(y)\,dv_*(y)\leqslant
C\cdot\int\limits_{\varepsilon_1}^{\varepsilon_2}\int\limits_{S(y_0,
r)}Q(y)\,d\mathcal{A}\,dr
\end{equation}
для кожних $0\leqslant \varepsilon_1<\varepsilon_2\leqslant
\delta(y_0).$
Покладемо
$r_*(y):=\min\{\varepsilon^{\,*}_1, r_0(y), \delta(y)\}.$
Покриємо множину $A_1$ кулями $B(y, r_*/4),$ $y\in A_1.$ Зауважимо,
що $|\gamma_1 |$ є компактною множиною в  ${\Bbb M}^n_*$ як
неперервний образ компакту $[1/2, 1]$ при відображенні $\gamma_1.$
Тоді за теоремою Гейне-Бореля-Лебега існує скінченне підпокриття
$\bigcup\limits_{i=1}^pB(y_i, r_*/4)$ множини $A_1.$ Іншими словами,
\begin{equation}\label{eq2A}
A_1\subset \bigcup\limits_{i=1}^pB(y_i,r_i/4)\,,\qquad 1\leqslant
p<\infty\,,
 \end{equation}
де $r_i$ позначає $r_*(y_i)$ при кожному $y_i\in \overline{D_*}.$

\medskip
Оскільки за припущенням випадку~1) $\gamma_m(t)\rightarrow
\partial D$ при $t\rightarrow c_m-0$ і $\Delta_m(t)\rightarrow
\partial D$ при $t\rightarrow d_m-0$ при всіх $m\in {\Bbb N},$
знайдуться послідовності точок $z^1_m\in |\gamma_m|$ and $z^2_m\in
|\Delta_m|$ такі, що
\begin{equation}\label{eq2D}
d(z^1_m,
\partial D)<1/m\,,\qquad  d(z^2_m,\partial D)<1/m\,,\quad
m=1,2,\ldots\,.
\end{equation}
Нехай $P_m$ -- частина кривої $\gamma_m$ в ${\Bbb M}^n,$ розташована
між точками $x_m$ і $z^1_m,$ і $Q_m$ -- частина кривої $\Delta_m$ в
${\Bbb M}^n,$ розташована між точками $x^{\,\prime}_m$ і $z^2_m.$ За
побудовою, $f(P_m)\subset A_1$ і $f(Q_m)\subset A_2.$ Покладемо
$\Gamma_m:=\Gamma(P_m, Q_m, D\setminus\{x_0\}).$ Надалі запис
$\Gamma_1>\Gamma_2$ означає, що кожна крива $\gamma_1\in \Gamma_1$
має підкриву $\gamma_2\in \Gamma_2.$ (Іншими словами, якщо
$\gamma_1\colon I\rightarrow{\Bbb M}^n,$ то $\gamma_2\colon
J\rightarrow{\Bbb M}^n,$ де $J\subset I$ і $\gamma_2(t)=\gamma_1(t)$
при $t\in J,$ а $I, J$ -- деякі відрізки, інтервали, або
півінтервали). Тоді з огляду на~(\ref{eq1J}) і (\ref{eq2A}), за
твердженням~\ref{pr2}
\begin{equation}\label{eq5}
\Gamma_m>\bigcup\limits_{i=1}^p\Gamma_{im}\,,
 \end{equation}
де $\Gamma_{im}:=\Gamma_{f_m}(y_i, r_i/4, r_i/2).$ Покладемо
$\widetilde{Q}(y)=\max\{Q(y), 1\}$ і
$\widetilde{q}_{y_i}(r)=\int\limits_{S(y_i,
r)}\widetilde{Q}(y)\,d\mathcal{A}.$ Тоді ми також отримаємо, що
$\widetilde{q}_{y_i}(r)\ne \infty$ для кожного $r\in [r_i/4,r_i/2].$
Покладемо
$$I_i=I_i(y_i,r_i/4,r_i/2)=\int\limits_{r_i/4}^{r_i/2}\
\frac{dr}{r\widetilde{q}_{y_i}^{\frac{1}{n-1}}(r)}\,.$$
Зауважимо, що $I\ne 0,$ оскільки $\widetilde{q}_{y_i}(r)\ne \infty$
при всіх $r\in [r_i/4,r_i/2].$ Окрім того, зауважимо, що
$I\ne\infty,$ бо для деякої сталої $C>0,$
$$I_i\leqslant C\cdot\log\frac{r_2}{r_1}<\infty\,,\quad i=1,2, \ldots, p\,.$$
Покладемо
$$\eta_i(r)=\begin{cases}
\frac{1}{I_ir\widetilde{q}_{y_i}^{\frac{1}{n-1}}(r)}\,,&
r\in [r_i/4,r_i/2]\,,\\
0,& r\not\in [r_i/4,r_i/2]\,.
\end{cases}$$
Зауважимо, що $\eta_i$ задовольняє
умову~$\int\limits_{r_i/4}^{r_i/2}\eta_i(r)\,dr=1,$ отже, її можна
підставити у праву частину нерівності~(\ref{eq2*A}) з відповідними
$r_1$ та $r_2.$ Ми будемо мати, що
\begin{equation}\label{eq7B}
M(\Gamma_{im})\leqslant \int\limits_{A(y_i, r_i/4, r_i/2)}
\widetilde{Q}(y)\,\eta^n_i(d_*(y, y_i))\,dv_*(y)\,.\end{equation}
Застосуємо оцінку~(\ref{eq7A}) у правій частині
співвідношення~(\ref{eq7B}). Ми отримаємо, що
$$\int\limits_{A(y_i, r_i/4, r_i/2)}
\widetilde{Q}(y)\,\eta^n_i(d_*(y, y_i))\,dv_*(y)\leqslant$$
\begin{equation}\label{eq7C}
\leqslant C_i \int\limits_{r_i/4}^{r_i/2}\int\limits_{S(y_i,
r)}Q(y)\eta^n_i(d_*(y, y_i))\,d\mathcal{A}\,dr\,=
\end{equation}$$=\frac{C_i}{I_i^n}\int\limits_{r_i/4}^{r_i/2}r^{n-1}
\widetilde{q}_{y_i}(r)\cdot
\frac{dr}{r^n\widetilde{q}^{\frac{n}{n-1}}_{y_i}(r)}=\frac{C_i}{I_i^{n-1}}\,,$$
де $C_i$ -- стала, що відповідає~$y_i$ у~(\ref{eq7A}). Тоді, з
огляду на~(\ref{eq7B}) і~(\ref{eq7C}), ми отримаємо, що
$M(\Gamma_{im})\leqslant \frac{C_i}{I_i^{n-1}},$
звідки з~(\ref{eq5})
\begin{equation}\label{eq7D}
M(\Gamma_m)\leqslant \sum\limits_{i=1}^pM(\Gamma_{im})\leqslant
\sum\limits_{i=1}^p\frac{C_i}{I_i^{n-1}}:=C_0\,, \quad
m=1,2,\ldots\,.
\end{equation}
Подальші міркування спираються на ``слабку плоскість'' ізольованих
точок межі області $D.$ Передусім, для будь-якої точки $z\in
\partial D,$ з огляду на співвідношення~(\ref{eq2D}) та нерівність
трикутника, при достатньо великих $m\in {\Bbb N}$ будемо мати:
$$d(x_m, z^1_m)\geqslant d(x_0, z)-d(z, z^1_m)-d(x_0, x_m)\geqslant$$
$$\geqslant d(x_0, z)-\frac{1}{m}-d(x_0, x_m)\geqslant \frac{1}{2}
d(x_0, z)\geqslant \frac{1}{2}d(x_0, \partial D)\,.$$
Міркуючи аналогічно, ми отримаємо, що
$d(x^{\,\prime}_m, z^2_m)\geqslant \frac{1}{2}d(x_0, \partial D).$
З огляду на отримані вище співвідношення,
$$\max\{d(P_m), d(Q(m)\}\geqslant \max \{d(x_m, z^1_m),
x^{\,\prime}_m, z^2_m)\}\geqslant \frac{1}{2}d(x_0, \partial D)$$
для достатньо великих $m\in {\Bbb N}.$

\medskip Нехай $U:=B(x_0, \delta_0/2),$ і нехай $V$ -- окіл цієї ж
точки $x_0,$ який відповідає твердженню~\ref{pr7} і
зауваженню~\ref{rem1}. Оскільки за припущенням $x_m, y_m\in
D\setminus \{x_0\},$ $m=1,2,\ldots,$ то знайдеться номер $m_0\in
{\Bbb N}$ такий, що $x_m, y_m\in V$ при всіх $m\geqslant m_0.$
Зауважимо, що для $m\geqslant m_0$
\begin{equation}\label{eq10B}
P_m\cap \partial V\ne\varnothing, \quad Q_m\cap
\partial V\ne\varnothing\,.
\end{equation}
Дійсно, $x_m\in P_m,$ $y_m\in Q_m,$ тому $P_m\cap V\ne\varnothing\ne
Q_m\cap V$ при $m\geqslant m_0.$ Крім того, ${\rm diam}\,V\leqslant
{\rm diam}\,U=\delta_0/2$ і, оскільки $d(P_m)\geqslant \delta_0>0$ і
$d(Q_m)\geqslant \delta_0>0$ при всіх $m\in {\Bbb N},$ то за
твердженням~\ref{pr2} є справедливими співвідношення~(\ref{eq10B}).
Аналогічно можна довести, що
\begin{equation}\label{eq11F}
P_m\cap \partial U\ne\varnothing, \quad Q_m\cap
\partial U\ne\varnothing\,.
\end{equation}
Тоді з огляду на твердження~\ref{pr7} і зауваження~\ref{rem1}
\begin{equation}\label{eq11G}
M(\Gamma_m)=M(\Gamma(P_m, Q_m, D\setminus \{x_0\}))>P\,,\qquad
m\geqslant m_0\,.
\end{equation}
Останнє співвідношення суперечить~(\ref{eq7D}), що завершує розгляд
випадку~\textbf{1)}.
\begin{figure}
\centerline{\includegraphics[scale=0.5]{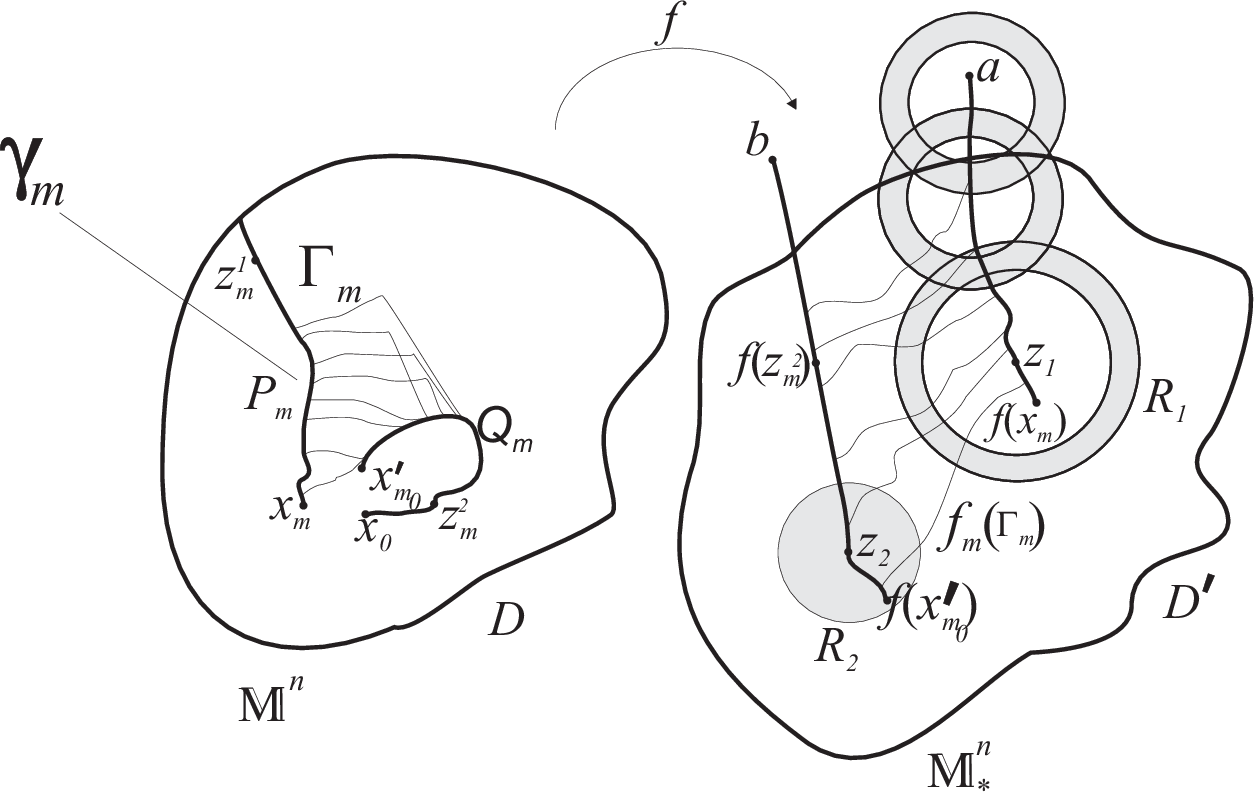}} \caption{До
доведення теореми~\ref{th1}. Випадок~2)}\label{figure2}
 \end{figure}
Розглянемо випадок \textbf{2)} при всіх $m\in {\Bbb N},$
$\gamma_m(t)\rightarrow
\partial D$ при $t\rightarrow c_m-0,$ але існує $m_0\in {\Bbb N}$
такий, що $\Delta_{m_0}(t)\rightarrow 0$ при $t\rightarrow
d_{m_0}-0,$ див. малюнок~\ref{figure2}.
Оскільки за припущенням $\gamma_m(t)\rightarrow
\partial D$ при $t\rightarrow c_m-0$ при всіх $m\in {\Bbb N}$ і
$\Delta_{m_0}(t)\rightarrow \{x_0\}$ при $t\rightarrow d_m-0,$
знайдуться послідовності точок $z^1_m\in |\gamma_m|$ and $z^2_m\in
|\Delta_{m_0}|$ такі, що
\begin{equation}\label{eq2E}
d(z^1_m,
\partial D)<1/m\,,\qquad  d(z^2_m, \{x_0\})<1/m\,,\quad
m=1,2,\ldots\,.
\end{equation}
Нехай $P_m$ -- частина кривої $\gamma_m$ в ${\Bbb M}^n,$ розташована
між точками $x_m$ і $z^1_m,$ і $Q_m$ -- частина кривої
$\Delta_{m_0}$ в ${\Bbb M}^n,$ розташована між точками
$x^{\,\prime}_{m_0}$ і $z^2_m.$ Ми бачимо, що ситуація у відомому
сенсі звелася до попереднього випадку~1), бо міркуючи тепер
аналогічно для $P_m$ і $Q_m$ і сім'ї кривих $\Gamma_m:=\Gamma(P_m,
Q_m, D\setminus\{x_0\}),$ як це було у попередньому випадку, ми
отримаємо з одного боку співвідношення~(\ref{eq7D}), а з іншого --
нерівність~(\ref{eq11G}), які суперечать одна одній.
Ми вже зауважили, що випадок~{\bf 3)} не вартий окремого розгляду,
оскільки відрізняється від випадку~{\bf 2)} лише позначеннями.

\begin{figure}
\centerline{\includegraphics[scale=0.5]{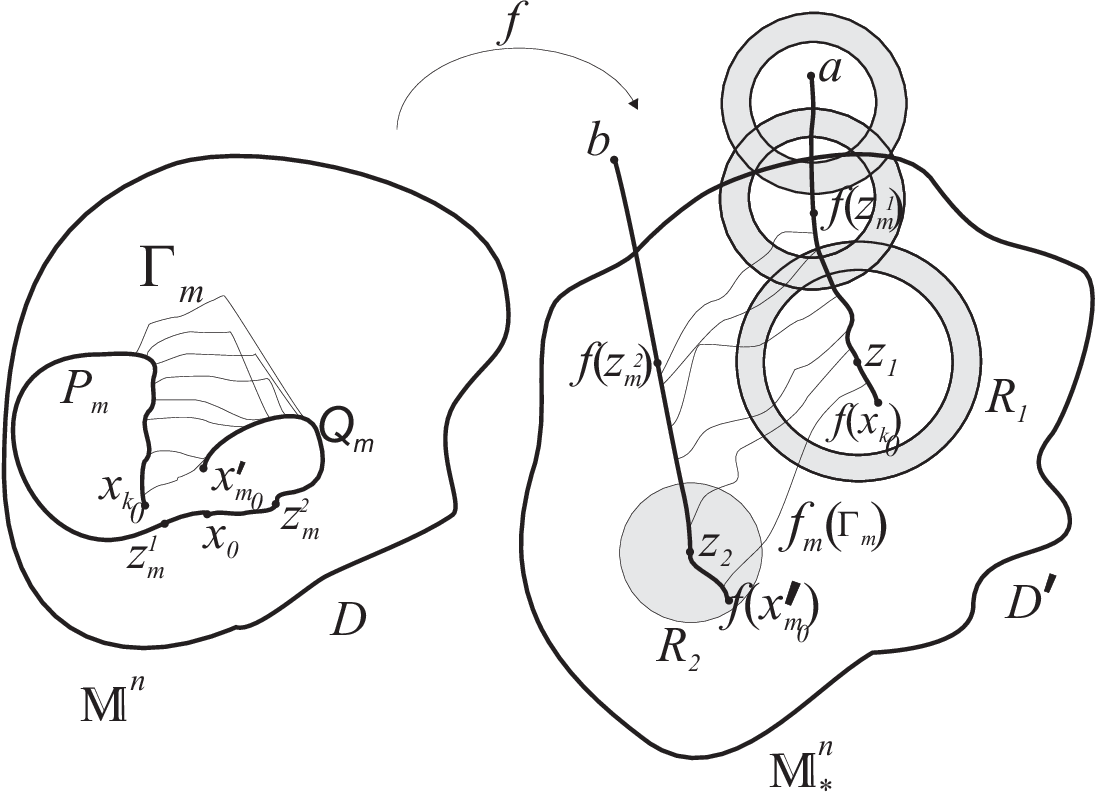}} \caption{До
доведення теореми~\ref{th1}. Випадок~4)}\label{figure3}
 \end{figure}

\medskip
Нарешті, розглянемо випадок~{\bf 4):} існують $k_0, m_0\in {\Bbb N}$
такі, що $\gamma_{k_0}(t)\rightarrow 0$ при $t\rightarrow c_{k_0}-0$
і $\Delta_{m_0}(t)\rightarrow 0$ при $t\rightarrow d_{m_0}-0,$ див.
малюнок~\ref{figure3}.
З умови випливає, що знайдуться послідовності точок $z^1_m\in
|\gamma_{k_0}|$ and $z^2_m\in |\Delta_{m_0}|,$ $m=1,2,\ldots $ такі,
що
\begin{equation}\label{eq2F} d(z^1_m,
\{x_0\})<1/m\,,\qquad  d(z^2_m, \{x_0\})<1/m\,,\quad m=1,2,\ldots\,.
\end{equation}
Нехай $P_m$ -- частина кривої $\gamma_{k_0}$ в ${\Bbb M}^n,$
розташована між точками $x_{k_0}$ і $z^1_m,$ а $Q_m$ -- частина
кривої $\Delta_{m_0}$ в ${\Bbb M}^n,$ розташована між точками
$x^{\,\prime}_{m_0}$ і $z^2_m.$ Як і вище, покладемо:
$\Gamma_m:=\Gamma(P_m, Q_m, D\setminus\{x_0\}).$ Міркуючи аналогічно
до випадку~1) для $P_m$ і $Q_m,$ ми знову отримаємо з одного боку
співвідношення~(\ref{eq7D}), а з іншого -- нерівність~(\ref{eq11G}),
які суперечать одна одній. Отже, наявність неперервного продовження
відображення $f$ у точку $x_0$ доведено.

\medskip
{\bf II. Відкритість та дискретність продовженого відображення.
Логарифмічна неперервність за Гельдером.} Нехай $U\subset {\Bbb
M}^n$ та $V\subset {\Bbb M}^n_*$ -- нормальні околи точок $x_0$ і
$f(x_0),$ відповідно. Нехай, крім того, $\varphi:U\rightarrow {\Bbb
R}^n$ та $\psi:V\rightarrow {\Bbb R}^n$ -- відповідні гомеоморфізми
цих околів у евклідовий простір. Можна вважати
$\varphi(x_0)=\psi(f(x_0))=0.$ Тоді $F(x):=\left(\psi\circ f\circ
\varphi^{\,-1}\right)(x)$ -- відображення між околами
$\varphi(U)\subset {\Bbb R}^n$ і $\psi(V)\subset {\Bbb R}^n$ початку
координат. За твердженням~\ref{pr4} $F$ є відкритим і дискретним,
тому таким є і відображення~$f.$

\medskip
Нехай $Q\in L^1(D^{\,\prime}).$ Тоді з огляду на
твердження~\ref{pr5} для кожної точки $y_0\in
\overline{D^{\,\prime}}$ знайдеться $r_0=r_0(y_0)>0$ таке, що
$q_{y_0}(r)<\infty$ при майже всіх $r\in (0, r_0);$ зокрема, за
пунктом ${\bf I}$ відображення $f$ має неперервне продовження у
точку $x_0,$ а за доведеним вище є відкритим і дискретним у $D.$
Крім того, за твердженням~\ref{pr6} відображення $F$ задовольняє
оцінку
\begin{equation}\label{eq3}
\biggl|\left(\psi\circ f\circ
\varphi^{\,-1}\right)(z)-\left(\psi\circ f\circ
\varphi^{\,-1}\right)(0)\biggr|= |F(z)|\leqslant\frac{C_n\cdot
(\Vert Q\Vert_1)^{1/n}}{\log^{1/n}\left(1+\frac{r_0}{|z|}\right)}
\end{equation}
для всіх $z\in \varphi(U).$ Слід зауважити, що з огляду означення
нормальних координат, для всяких $x, y\in U$ і $f(x), f(y)\in V$
$$C_1|\varphi(x)-\varphi(y)|\leqslant d(x, y)\leqslant C_2|\varphi(x)-\varphi(y)|\,,$$
\begin{equation}\label{eq4}
C_3|\psi(f(x))-\psi(f(y))|\leqslant d_*(f(x), f(y))\leqslant
C_4|\psi(f(x))-\psi(f(y))|\,,
\end{equation}
де $C_1, C_2, C_3, C_4$ -- деякі додатні сталі, близькі до~1 в
достатньо малих околах $U$ і $V.$ Покладемо $\varphi^{\,-1}(z)=x.$
Маючи на увазі, що $\varphi^{\,-1}(0)=x_0,$ з нерівності~(\ref{eq3})
ми отримаємо, що
\begin{equation}\label{eq5A}
\biggl|\left(\psi\circ f\right)(x)-\left(\psi\circ
f\right)(x_0)\biggr|\leqslant\frac{C_n\cdot (\Vert
Q\Vert_1)^{1/n}}{\log^{1/n}\left(1+\frac{r_0}{|\varphi(x)-\varphi(x_0)|}\right)}\,.
\end{equation}
З~(\ref{eq4}) та~(\ref{eq5A}) випливає, що
\begin{equation}\label{eq6}
\frac{1}{C_4}\cdot  d_*(f(x), f(x_0))\leqslant\frac{C_n\cdot (\Vert
Q\Vert_1)^{1/n}}{\log^{1/n}\left(1+\frac{r_0}{\frac{1}{C_1}\cdot
d(x, x_0)}\right)}\,.
\end{equation}
З огляду на правило Лопіталя
$\log^{1/n}\left(1+\frac{1}{nt}\right)\sim\log^{1/n}\left(1+\frac{1}{kt}\right)$
при $t\rightarrow+0$ для фіксованих $k, n> 0.$ Тоді з~(\ref{eq6})
випливає, що
$$d_*(f(x), f(x_0))\leqslant\frac{\widetilde{C_n}\cdot (\Vert
Q\Vert_1)^{1/n}}{\log^{1/n}\left(1+\frac{r_0}{d(x,
x_0)}\right)}\,,\quad x\in U\,.$$
Теорема повністю доведена.~$\Box$


КОНТАКТНА ІНФОРМАЦІЯ

\medskip
\noindent{{\bf Вікторія Сергіївна Десятка} \\
{\bf 1.} Житомирський державний університет ім.\ І.~Франко\\
кафедра математичного аналізу, вул. Велика Бердичівська, 40 \\
м.~Житомир, Україна, 10 008 \\
e-mail: victoriazehrer@gmail.com }

\medskip
\noindent{{\bf Євген Олександрович Севостьянов} \\
{\bf 1.} Житомирський державний університет ім.\ І.~Франко\\
кафедра математичного аналізу, вул. Велика Бердичівська, 40 \\
м.~Житомир, Україна, 10 008 \\
{\bf 2.} Інститут прикладної математики і механіки
НАН України, \\
вул.~Добровольського, 1 \\
м.~Слов'янськ, Україна, 84 100\\
e-mail: esevostyanov2009@gmail.com}


\begin{thebibliography}{99}

{\small

\bibitem{MRV$_1$} {\it Martio~O., Rickman~S., and V\"{a}is\"{a}l\"{a}~J.}
Definitions for quasiregular mappings // Ann. Acad. Sci. Fenn. Ser.
A1. -- 1969. -- \textbf{448.} -- P.~1--40.

\bibitem{Ri} {\it Rickman S.} Quasiregular mappings. -- Berlin: Springer-Verlag, 1993.

\bibitem{SevSkv$_1$} {\it Sevost'yanov~E.A., Skvortsov~S.A.} On mappings whose inverse
satisfy the Poletsky inequality // Ann. Acad. Scie. Fenn. Math. --
2020. -- \textbf{45.} -- P.~259--277.

\bibitem{Sev$_2$} {\it Sevost'yanov~E.A.}
On mappings with the inverse Poletsky inequality on Riemannian
manifolds // Acta Mathematica Hungarica. -- 2022. -- \textbf{167},
no.~2. -- P.~576--611.

\bibitem{Lee} {\it Lee J.\,M.} Riemannian Manifolds: An
Introduction to Curvature. -- New York: Springer, 1997.

\bibitem{MRSY} {\it Martio O., Ryazanov V., Srebro U. and Yakubov
E.} Moduli in Modern Mapping Theory. -- New York: Springer Science +
Business Media, LLC, 2009.

\bibitem{GSU} {\it Gol’dshtein~V., Sevost’yanov~E. and Ukhlov~A.}
On the boundary behavior of weak $(p; q)$-quasiconformal mappings //
Journal of Mathematical Sciences. -- 2023. -- \textbf{270.} --
P.~420-–427.

\bibitem{ARS} {\it Afanasieva~E.S., Ryazanov~V.I.,
Salimov~R.R.} On mappings in the Orlicz–Sobolev classes on
Riemannian manifolds // J. Math. Sci. -- 2012. -- \textbf{181},
no.~1. -- P.~1--17.

\bibitem{SM} {\it Sevost'yanov~E.A. and Markysh~A.A.}
On Sokhotski--Casorati--Weierstrass theorem on metric spaces //
Complex Variables and Elliptic Equations. -- 2019. -- \textbf{64},
no.~12. -- P.~1973--1993.

\bibitem{Va} {\it V\"{a}is\"{a}l\"{a} J.} Lectures on $n$-Dimensional Quasiconformal
Mappings. --  Lecture Notes in Math. 229, Berlin etc.:
Springer--Verlag, 1971.

\bibitem{Fe} {\it Federer~H.:} Geometric Measure Theory. --
Berlin etc.: Springer, 1969.

\bibitem{Ku} {\it Kuratowski K.} Topology, v. 2. --
New York--London: Academic Press, 1968.

\bibitem{MRV$_3$} {\it Martio~O., Rickman~S., and
V\"{a}is\"{a}l\"{a}~J.} Topological and metric properties of
quasiregular mappings // Ann. Acad. Sci. Fenn. Ser. A1. -- 1971. --
\textbf{488}. -- P.~1--31.

\bibitem{Skv} {\it Скворцов С.О.} Локальна поведiнка вiдображень метричних просторiв з розгалуженням
// Український математичний вiсник. -- 2020. -- \textbf{17}, no.~4. -- P.~574--593;
transl. ``Local behavior of mappings of metric spaces with
branching'' in Journal of Mathematical Sciences. -- 2021. --
\textbf{254}, No. 3. -- P.~425--574.

\bibitem{HW} {\it Hurewicz~W., Wallman~H.} Dimension theory. -- Princeton: Princeton Univ. Press,
1948.

\bibitem{TY} {\it Titus~C.J. and Young~G.S.} The extension of interiority with some
applications // Trans. Amer. Math. Soc. -- 1962. -- \textbf{103.} --
P.~329--340.

\bibitem{IS} {\it Ilyutko D. and Sevost'yanov E.}
On local properties of one class of mappings on Riemannian manifolds
// Journal of Mathematical Sciences. -- 2015. -- \textbf{211}, no.
5. -- P. 660--667.

\bibitem{SSD} {\it Sevost'yanov~E.A., Skvortsov~S.O., Dovhopiatyi~O.P.}
On nonhomeomorphic mappings with the inverse Poletsky inequality //
Journal of Mathematical Sciences. -- 2021. -- \textbf{252}, no.~4.
-- P.~541--557.

\bibitem{IS$_2$} {\it Ilyutko D. and Sevost'yanov E.}
On the equicontinuity of families of inverse mappings of Riemannian
manifolds // Journal of Mathematical Sciences. -- 2020. --
\textbf{246}, no. 5. -- P. 664--670.

\bibitem{SevSkv$_2$} {\it Севостьянов~Е.А., Скворцов~С.А.} О локальном поведении
одного класса обратных отображений // Укр. мат. вестник. -- 2018. --
Т.~\textbf{15}, №~3. -- С.~399--417; translation ''On the local
behavior of a class of inverse mappings'' in J. Math. Sci. -- 2019.
-- V.~\textbf{241,} no.~1. -- P.~77--89.

\bibitem{Sev$_3$} {\it Sevost'yanov~E.A.} Isolated singularities of mappings with the inverse Poletsky
inequality // Mat. Stud. -- 2021. -- \textbf{55}, no.~2. --
P.~132--136. }


\end{thebibliography}
\end{document}